\documentclass[letterpaper, 10 pt, conference]{ieeeconf}  

\IEEEoverridecommandlockouts                              
\usepackage{amsfonts}
\usepackage{graphicx,xcolor}
\usepackage{amsmath}
\usepackage{amssymb}
\usepackage{dsfont}

\usepackage{ifthen}
\usepackage{color}
\usepackage{cite}

\usepackage{mathtools}
\usepackage{booktabs}
\usepackage{url}
\usepackage{hyperref}
\usepackage{caption}
\usepackage{subfig}
\usepackage{float}
\usepackage{multirow}

\def\scr#1{{\cal #1}}
\newcommand{\R}{{\rm I\!R}}
\def\eq#1{\begin{equation}#1\end{equation}}
\def\rep#1{(\ref{#1})}
\newcommand{\bbb}{\mathbb}
\newtheorem{theorem}{Theorem}
\newtheorem{lemma}{Lemma}

\newtheorem{example}{Example}
\newtheorem{proposition}{Proposition}
\newtheorem{corollary}{Corollary}

\newcommand{\dfb}{\stackrel{\Delta}{=}}
\def\qed{ \rule{.08in}{.08in}}

\newcommand{\1}{\mathbf{1}}

\newcommand{\y}{\mathbf{y}}

\title{\LARGE \bf Reaching a Consensus with Limited Information
}

\author{Jingxuan Zhu \hspace{.3in} Yixuan Lin \hspace{.3in} Ji Liu \hspace{.3in} A. Stephen Morse
\thanks{J.~Zhu and Y.~Lin are with the Department of Applied Mathematics and Statistics at Stony Brook University (\texttt{\{jingxuan.zhu,yixuan.lin.1\}@stonybrook.edu}).
J. Liu is with the Department of Electrical and Computer Engineering at Stony Brook University
(\texttt{ji.liu@stonybrook.edu}).
A.S. Morse is with the Department of Electrical Engineering at Yale University
(\texttt{as.morse@yale.edu}).
}
}

\begin{document}

\maketitle
\thispagestyle{empty}
\pagestyle{empty}


\begin{abstract}
In its simplest form the well known consensus problem for a networked family of autonomous agents is to devise a set of protocols or update rules, one for each agent, which can enable
all of the agents to adjust or tune their ``agreement variable'' to the same value by utilizing real-time information obtained from their ``neighbors'' within the network. The aim of this paper is to study the problem of achieving a consensus in the face of limited
information transfer between agents. By this it is meant that instead of each agent receiving an agreement variable or real-valued state vector from each of its neighbors, it receives a linear function of each state instead. The specific problem of interest is formulated and provably correct algorithms are developed for a number of
special cases of the problem.
\end{abstract}

\vspace{.1in}

\section{Introduction}


In its simplest form the well known  consensus problem
\cite{lamport} for a networked family  of autonomous agents is
to devise a set of protocols or update rules, one for each
agent, which can enable all of the  agents to adjust  or tune
their ``agreement variable''  to the same value by utilizing
real-time
 information obtained from their ``neighbors'' within the network. The  consensus problem is  one
  of the most fundamental problems in the area of distributed computation and control.
 Consensus algorithms can be found as components of
 a  large variety of  more specialized algorithms  in the area of  distributed  computation and control
 such as distributed algorithms   for solving linear algebraic equations \cite{tacle}, distributed
 optimization problems \cite{nedic2009distributed}, distributed estimation problems \cite{kalman}, and even some distributed
 control problems \cite{overview}.

 There are a great many variations of the consensus problem. For example, the agreement variables
 could be restricted to be  real-valued  vectors or alternatively  integer-valued vectors \cite{basar}. The updating
  of agreement variables could be executed either  synchronously or asynchronously \cite{async}. The
   topology of the  network could be fixed or changing with time \cite{luc}. There could be malicious
    agents attempting to prevent consensus \cite{fischer}. There could be communication delays \cite{Ts84}
    or bit-rate constraints \cite{quant}. The target value of the agreement variables could be unconstrained  or it
    could  be some specified function of the initial values of the agents' agreement variables as for
     example in    distributed averaging \cite{Boyd2004} or gossiping \cite{deterministic}.  Some versions of the problem
      such as when agreement variables  take values in a finite set, defy deterministic solutions
       \cite{basar} whereas other  versions of the problem do not.

The aim of this paper is to  study  the problem of achieving a consensus  in the face of  limited
information transfer between   agents. The problem setup is as follows.  We consider a group of
$m>1$ autonomous agents labeled $1$ to $m$. Each agent  $i$ has a set of neighbors
 from whom agent $i$  can receive information; the set of labels of agent $i$'s neighbors
 (excluding itself), denoted by    $\scr{N}_i\subset\mathbf{m}\dfb \{1,2,\ldots,m\}$,\footnote{We use $\scr A\subset \scr B$ to denote that $\scr A$ is a subset of $\scr B$.}  is part of
  the problem formulation. The neighbor sets $\scr{N}_i,\;i\in\mathbf{m}$,  determine an $m$-vertex
   directed graph $\mathbb{N}$    defined so that there is an arc (or a directed edge) from vertex $j $ to vertex $i$ just
    in case agent $j$ is a neighbor of agent $i$.
Each agent $i$ has an agreement variable  or state $x_i\in\R^n$ which it can adjust synchronously
at times $t\in\{0,1,2,\ldots\}$.
 At time $t$,   agent $i$ receives from each neighbor $j\in\scr{N}_i$ a signal $s_{ji}(t) = C_{ji}x_j(t)$
 where $C_{ji}$ is a fixed real-valued matrix. Associating each arc $(j,i)$ in $\bbb N$ with matrix $C_{ji}$ leads to a matrix-valued weighted neighbor graph $\bar{\bbb N}$. It is assumed
  that for each $j\in\scr{N}_i$, $i\in\mathbf{m}$, both agents $i$ and $j$ know  $C_{ji}$. There are no
  priori constraints on $C_{ji}$. Some could, for example, be matrices with less rows then columns
  in which cases the information transferred by each such  corresponding signal $s_{ji}(t) = C_{ji}x_j(t)$
  would be insufficient to  determine $x_j(t)$. 
  In this sense the information
 agent $i$ receives from neighbor $j$ at time $t$ is limited to only a ``part of''  $x_j(t)$. Given this setup,
the  consensus problem of  interest is to devise  update rules using the $s_{ji}(t)$, one for each agent,
  which if possible will cause all $m$ agents' states $x_i$, $i\in\mathbf{m}$,   to converge to the same value in the limit
   as $t\rightarrow\infty$.

\vspace{.1in}

\section{Well-configured Systems} \label{sec:well system}

Consider the multi-agent system just described. 
We say that the $m$ agents are in {\em local agreement} with specific states $x_i$, $i\in\mathbf{m}$, if 
$C_{ji}x_i=C_{ji}x_j$ for all $i\in\mathbf{m}$ and $j\in\scr N_i$. 
We say that the $m$ agents have
reached a {\em consensus} with specific states $x_i$, $i\in\mathbf{m}$, if $x_i=x_j$ for all $i,j\in\mathbf{m}$.
A weighted neighbor graph $\bar{\bbb N}$ is called {\em well-configured} if local agreement implies consensus.

A well-configured weighted neighbor graph $\bar{\bbb N}$ has the following equivalent mathematical description. For each vertex $i$ in $\bbb N$, let $d_i$ denote the number of neighbors of agent $i$. Then $d=\sum_{i=1}^m d_i$ equals the total number of directed edges in $\scr E$.
Let $k_{i1}, \ldots, k_{id_i}$ be an arbitrary ordering of the labels in $\mathcal{N}_i$. Label all the $d$ arcs from 1 to $d$ according to the sequence $k_{11}, \ldots, k_{1d_1},\ldots,k_{m1},\ldots, k_{md_m}$. 
Define the corresponding {\em incidence matrix} $J$ as an $m\times d$ matrix in which column $k$ has exactly one $1$ in row $i$ and exactly one $-1$ is row $j$ if the $k$th arc in $\bbb N$ is $(j, i)$. 
For any finite set of matrices $\{M_1,M_2,\ldots,M_k\}$, we use ${\rm blockdiag}\{M_1,M_2,\ldots,M_k\}$ to denote the block diagonal matrix whose $i$th diagonal block is $M_i$. Define
\begin{align*}
    C &= {\rm blockdiag}\big\{ C_{k_{11},1}, \ldots, C_{k_{1d_1},1}, \ldots, \\
    & \;\;\;\;\;\;\;\;\;\;\;\;\;\;\;\;\;\;\;\;\;\;\;\;\;\;\;\;\;\;\;\;\;\;\;\;\;\;\;\;\;\;\;\;\;C_{k_{m 1},m}, \ldots, C_{k_{m d_m},m}\big\}. 
\end{align*}
Let $\bar J = J \otimes I_n$  and $\bar I = \1_m\otimes I_n$, 
where $\otimes$ denotes the Kronecker product, $I_n$ denotes the $n\times n$ identity matrix, and $\1_m$ denotes the $m$-dimensional column vector whose entries all equal 1. Then it is not hard to verify that
a weighted neighbor graph $\bar{\bbb N}$ is well-configured if and only if 
\eq{{\rm kernel}\;C\bar J'={\rm span}\;\bar I.\label{eq:well}}
In the case when $\bbb N$ is weakly connected\footnote{A directed graph is weakly connected if there is an undirected path between
each pair of distinct vertices.}, $\text{kernel}\;\bar{J}' = \text{span}\;\bar{I}$ \cite[Theorem 8.3.1]{algebraicgraph}; then \eqref{eq:well} will be true if and only if 
\eq{\text{span}\;\bar{J}'\cap\text{kernel}\;C = 0.\label{jiliu}}
It is worth emphasizing that $C$ and $J$ are defined according to the same ordering of the arcs in $\bar{\bbb N}$, and the necessary and sufficient condition \eqref{eq:well} or \rep{jiliu} is independent of the ordering. 


With the above in mind, the following two questions arise. First, what are the necessary and/or sufficient conditions on $\bbb N$ for which there exist $C_{ji}$ matrices so that $\bar{\bbb N}$
is well-configured?
Second, if $\bar{\bbb N}$ is well-configured, how one can construct a recursive distributed algorithm for each agent which will drive the system
from arbitrary start states to local agreement and thus to a consensus? These are precisely what we consider in this paper.

\vspace{.1in}

\section{System Design}

The goal of this section is to derive graph-theoretic conditions on which a multi-agent system can be well-configured. 

As described, for any pair of neighboring agents, say agent $i$ and its neighbor $j$, agent $j$ only sends $C_{ji}x_j$ to agent $i$ so that the transmitted vector size may be reduced and $x_j$ may not be identified. Thus it is sometimes desirable that $\scr K_{ji}\neq 0$, where $\scr K_{ji}$ denotes the kernel of $C_{ji}$; otherwise, $x_j$ can be uniquely determined from $C_{ji}x_j$. Also, if $\scr K_{ji}\neq 0$, the size of $C_{ji}x_j$ will be no smaller than that of $x_j$.

A directed graph $\bbb G$ is called rooted if it contains a directed spanning tree of $\bbb G$, and called strongly connected if there is a directed path between each pair of distinct vertices. Every strongly connected graph is rooted, but not vice versa.

First, it is easy to see that if $\bbb N$ is not rooted, a consensus cannot be guaranteed for arbitrary initial values. We next consider some examples of rooted graphs. 


\subsection{Rooted Graphs}

If $\bbb N$ is rooted, $\bar{\bbb N}$ cannot be always well-configured with all $\scr K_{ji}\neq 0$, as shown in the following lemma for path graphs. 

\vspace{.05in}

\begin{lemma}\label{lm:path}
If $\bbb N$ is a directed path, then $\bar{\bbb N}$ can be well-configured only if all $\scr K_{ji}=0$. 
\end{lemma}

\vspace{.05in}

{\bf Proof of Lemma \ref{lm:path}:}
For a directed path $x_1, x_2, \ldots, x_m,$ we have local agreement $C_i(x_i - x_{i+1}) = 0,$ for $i=1, \ldots, m-1.$ Suppose to the contrary that there exists an $i$ such that $\scr K_i\neq 0$, then there exists a nonzero $x$ such that $C_ix = 0.$ Let 
\begin{align*}
    x_j &= x_1,\,\,\, j = 1,\ldots, i\\
    x_j &= x_1 + x, \,\,\, j = i+1,\ldots,m.
\end{align*}
We have $C_i(x_i - x_{i+1}) = 0$ for $i=1, \ldots, m-1,$ while $x_i$ do not reach a consensus.
\hfill$\qed$

\vspace{.05in}

The following example shows that there exists a rooted graph which can be well-configured with all $\scr K_{ji}\neq 0$. 

\vspace{.05in}

\begin{example}
Consider a three-agent network with arcs $1\rightarrow2, 2\rightarrow1, 3\rightarrow1, 3\rightarrow2$. Then local agreement equations are
\begin{equation}\label{eq:example1}
\begin{split}
    C_{12}(x_1-x_2)&=0 \\
    C_{21}(x_2-x_1)&=0 \\
    C_{31}(x_3-x_1)&=0 \\
    C_{32}(x_3-x_2)&=0 
\end{split}
\end{equation}
Note that the existence $x_i$ satisfying the four equalities above imply that there are vectors $p_1, p_2, p_3$, namely 
$p_1 = x_1 - x_2$, $p_2 = x_3 - x_1$, $p_3 = x_2 - x_3,$ 
such that
\begin{equation}\label{eq:examp1_1}
    \begin{split}
    p_1 + p_2 + p_3 &= 0 \\
    p_1 &\in \scr K_{12}\cap\scr K_{21} \\
    p_2 &\in \scr K_{31} \\
    p_3 &\in \scr K_{32}
    \end{split}
\end{equation}
Conversely, for any set of vectors $p_1, p_2, p_3$ satisfying \eqref{eq:examp1_1}, there are vectors $x_i$, namely $x_1 =
p_1$, $x_2 = 0$, $x_3 = p_3$ for which the four equalities in \eqref{eq:example1} hold.
Note that there will exist $p_1, p_2, p_3$ for which \eqref{eq:examp1_1} holds if and only if
\begin{align*}
    p_1 &\in \scr K_{12}\cap\scr K_{21} \cap (\scr K_{31}+\scr K_{32}) \\
    p_2 &\in \scr K_{31}\cap (\scr K_{12} \cap \scr K_{21}+\scr K_{32}) \\
    p_3 &\in \scr K_{32}\cap (\scr K_{12} \cap \scr K_{21}\cap\scr K_{31})
\end{align*}
Thus the conditions for the $p_i$ to all equal zero are
\begin{align*}
    \scr K_{12}\cap\scr K_{21} \cap (\scr K_{31}+\scr K_{32}) & = 0\\
    \scr K_{31}\cap (\scr K_{12} \cap \scr K_{21}+\scr K_{32}) & = 0\\
    \scr K_{32}\cap (\scr K_{12} \cap \scr K_{21}\cap\scr K_{31}) & = 0
\end{align*}
which are the conditions for the three subspaces $\scr K_{12}\cap\scr K_{21}$, $\scr K_{31}$ and $\scr K_{32}$ to be independent. Thus the weighted neighbor graph of interest is well-configured just in case the three subspaces are independent, and do not necessarily have to equal 0.
\hfill$\Box$
\end{example}

\vspace{.05in}

It turns out that well-configuration characterization of rooted graphs is quite complicated. We thus leave it as a future direction and focus on strongly connected graphs in the next subsection. 

\subsection{Strongly Connected Graphs}

Strong connectedness itself cannot guarantee well-configuration. 
To state our sufficient condition for well-configuration, we need the following concept from graph theory \cite{ear}. 

\begin{figure*}[!ht]
\centering
\includegraphics[width=1.0\textwidth]{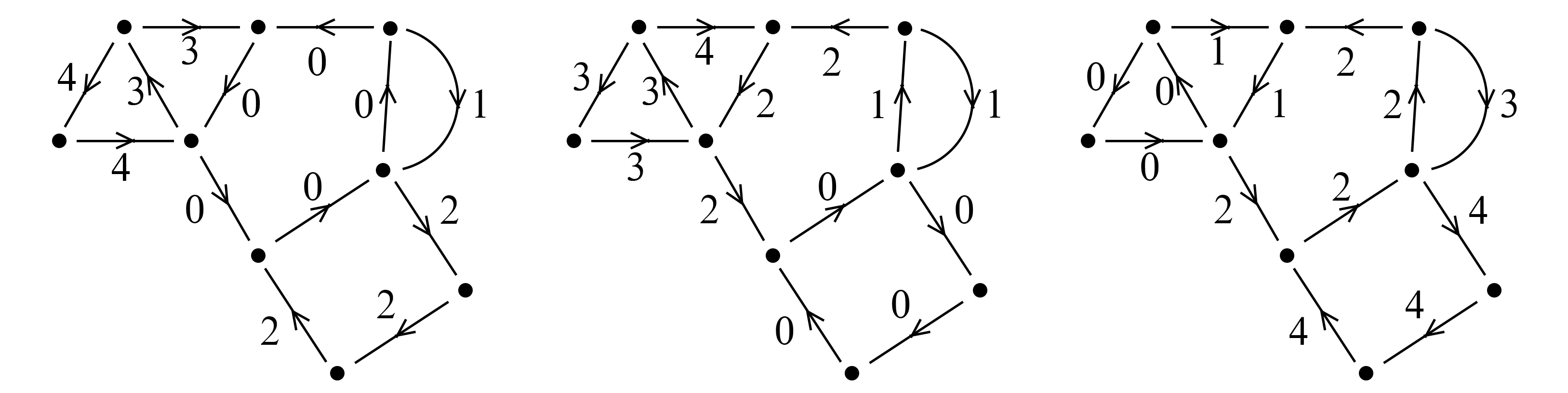} 
\caption{An example of different ear decompositions of a strongly connected graph. The number associated to each arc represents the index of the ear which the arc belongs to in an ear decomposition. }
\label{fig:reduction}
\end{figure*}

An {\em ear decomposition} of a directed graph without self-arcs\footnote{The definition can be extended to more general directed multigraphs with self-arcs \cite{ear}.} $\bbb G = (\scr V, \scr E)$ with at least two vertices is a sequence of subgraphs of $\bbb G$, denoted $\{\bbb E_0,\bbb E_1,\ldots,\bbb E_p\}$, in which $\bbb E_0$ is a directed cycle, and each $\bbb E_i$, $i\in\mathbf{p}$, is a directed path or a directed cycle with the following properties:
\begin{enumerate}
   \item $\{\bbb E_0,\bbb E_1,\ldots,\bbb E_p\}$ form an arc partition of $\bbb G$, i.e., 
   $\bbb E_i$ and $\bbb E_j$ are arc disjoint if $i\neq j$, and $\bigcup_{k=0}^{p}\bbb E_k = \bbb G$;
   \item For each $i\in\mathbf{p}$, if $\bbb E_i$ is a directed cycle, then it has precisely one vertex in common with $\bigcup_{k=0}^{i-1}\bbb E_k$; if $\bbb E_i$ is a directed path, then its two end-vertices are the only two vertices in common with $\bigcup_{k=0}^{i-1}\bbb E_k$.
\end{enumerate}
Each of $\bbb E_0,\bbb E_1,\ldots,\bbb E_p$ is called an {\em ear} of the decomposition.
Not all directed graphs admit an ear decomposition. It has been proved that a directed graph has an ear decomposition if and only if it is strongly connected \cite[Theorem 7.2.2]{ear}. It is also known that there exists a linear algorithm to find one ear decomposition of a strongly connected graph \cite[Corollary 7.2.5]{ear}. A strongly connected graph may admit multiple ear decompositions, and apparently, the number of all possible different ear decompositions of a strongly connected graph is finite. It turns out that every ear decomposition of a strongly connected graph with $m$ vertices and $e$ arcs has $e-m+1$ ears \cite[Corollary 7.2.3]{ear}.
To help understand the concept, an illustrative example is provided in Figure \ref{fig:reduction}.










Two subspaces $\scr S_1$ and $\scr S_2$ of $\R^n$ are independent if their intersection is the zero subspace, i.e., if $\scr S_1\bigcap \scr S_2 = 0$.
A finite family of subspaces $\{\scr S_1,\scr S_2,\ldots, \scr S_p\}$  is independent if $$\scr S_i\bigcap \Big(\sum_{j\neq i} \scr S_j\Big)=0, \;\;\;\;\;i\in\mathbf{p}.$$


\begin{theorem}\label{thm:steve}
Suppose that $\bbb N$ is strongly connected and let $\mathrm{D}$ be an ear decomposition of $\bbb N$. If for each ear $\bbb E\in\mathrm{D}$, $\{\scr K_{ji} : (j,i)\in\bbb E\}$ is an independent family, then $\bar{\bbb N}$ is well-configured. 
\end{theorem}

\vspace{.05in}

To prove the theorem, we first study directed cycles and paths since they are basic components in ear decompositions. 


To simplify notation, we label the vertices of an $m$-vertex directed cycle as $1\rightarrow 2\rightarrow\cdots\rightarrow m\rightarrow1$. Suppose that $C_1,C_2,\ldots,C_m$ are
given matrices, each with $n$ columns. Suppose that for each $i\in\mathbf{m}$, agent $i$ receives $C_ix_{i-1}$ from agent $i-1$, where it is understood that agent $0$ and agent $m$ are one and the same, and that $x_0\dfb x_m$. Thus for this $\bar{\bbb N}$ to be well-configured means that the relations
\eq{C_ix_i = C_ix_{i-1}, \;\;\;\;\; i\in\mathbf{m},\label{eq:cycle1}}
must imply that $x_i = x_{i-1}$, $i\in\mathbf{m}$.
Let $\scr K_i$ denote the kernel of $C_i$ for all $i\in \mathbf{m}$.


\begin{lemma}\label{lem:cycle}
If $\bbb N$ is an $m$-vertex directed cycle, then $\bar{\bbb N}$ is well-configured by matrices $C_i$, $i\in\mathbf{m}$, if and only if $\{\scr K_1,\scr K_2, \ldots, \scr K_m\}$ is an independent family.
\end{lemma}

\vspace{.05in}

{\bf Proof of Lemma \ref{lem:cycle}:}
Since $\bbb N$ is a cycle, $m\ge 2$.
We first prove the sufficiency. Let $\{\scr K_i,\; i\in\mathbf{m}\}$ be an independent family. Suppose to the contrary that $\bar{\bbb N}$ is not well-configured. Then there must exist a non-consensus set $\{x_1,\ldots,x_m\}$ which satisfies \eqref{eq:cycle1}. Let $y_i = x_i - x_{i+1}$ for each $i \in \{1,\ldots, m-1\}$ and $y_m = x_m -x_1$. Then at least one of $y_1,\ldots,y_m$ is nonzero. Since $\sum_{i = 1}^m y_i = 0$, at least two of $y_1,\ldots,y_m$ are nonzero. Let $\scr A = \{i\in\mathbf{m}: y_i\neq 0\}$. Then $|\scr A|\ge 2$ and $\sum_{i\in\scr A} y_i = 0$. Since each $y_i\in\scr K_i$ and $\{\scr K_i, \; i\in\scr A\}$ is an independent family, $y_i$, $i\in\scr A$, are linear independent, which is contradictory to $\sum_{i\in\scr A} y_i = 0$. 

We next prove the necessity. Let $\{C_1,\ldots,C_m\}$ be any set of matrices which make $\bar{\bbb N}$ well-configured. Suppose to the contrary that $\{\scr K_1,\ldots,\scr K_m\}$ is not an independent family, which implies that there exists an index $p$ such that
$\scr K_p \cap (\sum_{i\neq p} \scr K_i)\neq 0$.
Then there exist $k_i\in\scr K_i$, $i\in\mathbf{m}$, such that 
$k_p=\sum_{i\neq p} k_i$,
which is nonzero. For any $x_1$, let $x_{i+1}=x_i+k_i$ for each $i\in\{1,\cdots, p-1,p+1,\ldots,m-1\}$ and $x_{p+1}=x_p-k_p$. It is easy to see check that such a set of non-consensus vectors $x_1,\ldots,x_m$ satisfy \eqref{eq:cycle1}. But this is impossible as $\bar{\bbb N}$ is well-configured.
\hfill$\qed$

\vspace{.05in}

It is easy to see that $n$ is the maximum possible number of subspaces in an independent family of nonzero subspaces of~$\R^n$. We thus have the following immediate consequence of Lemma~\ref{lem:cycle}.

\vspace{.05in}

\begin{corollary}\label{prop:cycle}
If $\bbb N$ is an $m$-vertex directed cycle, then $\bar{\bbb N}$ can be well-configured with all $\scr K_{i}\neq 0$, $i\in\mathbf{m}$, if and only if 
$m\le n$.
\end{corollary}

\vspace{.04in}

More can be said. 

\vspace{.04in}

\begin{lemma}\label{lem:reducedcycle}
Let $\bbb N$ be an $m$-vertex directed cycle with edge set $\scr E_{\bbb N}$. Let $\scr E$ be a subset of $\scr E_{\bbb N}$ defined as $\scr E=\{(i,j)\in\scr E_{\bbb N} : x_i=x_j\}$.
Then $\bar{\bbb N}$ is well-configured by matrices $C_i$, $i\in\mathbf{m}$, if and only if $\{\scr K_i : i\in\mathbf{m}, (i-1,i)\notin \scr E\}$ is an independent family.
\end{lemma}

\vspace{.05in}

{\bf Proof of Lemma \ref{lem:reducedcycle}:}
The case of $m-|\scr E|=0$ is trivial. 
We claim that $m-|\scr E|\neq 1$. To see this, suppose to the contrary that $m-|\scr E|=1$. Then the edges in $\scr E$ forms a directed spanning path of $\bbb N$, which guarantees that all $m$ agents reach a consensus. This implies that all the edges of $\bbb N$ belong to $\scr E$, which is impossible. Thus we focus on $m-|\scr E|\ge 2$ in the remaining proof.

Each vertex $i$ in directed cycle $\bbb N$ has a unique outgoing neighbor, denoted as $v[i]$. Let $\scr V$ be the vertex subset defined as $\scr V=\{i\in\mathbf{m}  :  (i,v[i])\notin\scr E\}$. Then $|\scr V|=m-|\scr E|$. 
Relabel the vertices in $\scr V$ as $v_1,\ldots,v_p$, $p=m-|\scr E|\ge 2$, along with the same direction as the directed cycle. It is not hard to verify that 
\begin{equation*}
    \begin{split}
        C_{v_1}(x_{v_1} - x_{v_2}) &= 0\\
    &\vdots\\
    C_{v_{p-1}}(x_{v_{p-1}} - x_{v_p}) &= 0\\
    C_{v_p}(x_{v_p} - x_{v_1}) &= 0
    \end{split}
\end{equation*}
which are mathematically equivalent to \eqref{eq:cycle1} with $m$ being replaced by $p$. Thus the above equations are equivalent to local agreements of an $p$-vertex directed cycle. 
From Lemma~\ref{lem:cycle}, $\bar{\bbb N}$ can be well-configured with all $\scr K_{v_i}\neq 0$, $i\in\mathbf{p}$, if and only if $p\le n$.
\hfill$\qed$

\vspace{.05in}

Lemma \ref{lem:reducedcycle} immediately implies the following result. 

\vspace{.05in}

\begin{corollary}\label{prop:reducedcycle}
Let $\bbb N$ be an $m$-vertex directed cycle with edge set $\scr E_{\bbb N}$. Let $\scr E$ be a subset of $\scr E_{\bbb N}$ defined as $\scr E=\{(i,j)\in\scr E_{\bbb N} : x_i=x_j\}$.  
Then $\bar{\bbb N}$ can be well-configured with all $\scr K_i\neq 0$ if and only if $m-|\scr E|\le n$. 
\end{corollary}

\vspace{.05in}

The above results can be directly applied to the following special case of path graphs. 

To simplify notation, we label the vertices of an $m$-vertex directed path as $1\rightarrow 2\rightarrow\cdots\rightarrow m$. Suppose that $C_2,\ldots,C_m$ are
given matrices, each with $n$ columns. Suppose that for each $i\in\mathbf{m}$, agent $i$ receives $C_ix_{i-1}$ from agent $i-1$. Thus for this $\bar{\bbb N}$ to be well-configured means that the relations $C_ix_i = C_ix_{i-1}$, $i\in\{2,\ldots,m\}$,
must imply that 
$x_i = x_{i-1}$, $i\in\{2,\ldots,m\}$.
Adding the arc $(m,1)$ to the above path and imposing $x_1=x_m$ will lead to a special case satisfying the condition in Lemma \ref{lem:reducedcycle} and Corollary \ref{prop:reducedcycle}, which immediately implies the following result.

\vspace{.05in}

\begin{corollary}\label{coro:path}
If $\bbb N$ is an $m$-vertex directed path with $x_1=x_m$, then $\bar{\bbb N}$ is well-configured by matrices $C_i$, $i\in\{2,\ldots,m\}$, if and only if $\{\scr K_2, \ldots, \scr K_m\}$ is an independent family, and thus $\bar{\bbb N}$ can be well-configured with all $\scr K_{i}\neq 0$, $i\in\{2,\ldots,m\}$, if and only if 
$m-1\le n$.
\end{corollary}

\vspace{.05in}

Compared with Lemma \ref{lm:path}, it is worth emphasizing that assuming $x_1=x_m$ significantly changes the condition for well-configuration of path graphs.

We are now in a position to prove Theorem \ref{thm:steve}.

\vspace{.05in}

{\bf Proof of Theorem \ref{thm:steve}:}
Let $\mathrm{D}=\{\bbb E_0,\bbb E_1,\ldots,\bbb E_p\}$ be the given ear decomposition of $\bbb N$. We claim that for each $i\in\{0,1,\ldots,p\}$, $\bigcup_{k=0}^{i}\bbb E_k$ is well-configured. The claim will be proved by induction on the index $i$.

By definition, $\bbb E_0$ is a directed cycle. From Lemma~\ref{lem:cycle}, $\bbb E_0$ is well-configured. 
Now suppose that the claim holds for all $i$ in the range $0\le i \le j$, where $j$ is a nonnegative integer
smaller than $p$. Consider ear $\bbb E_{i+1}$, which is either a directed cycle or a directed path. We treat  these two cases separately. If $\bbb E_{i+1}$ is a directed cycle, using the preceding argument, it is well-configured. Since ear $\bbb E_{i+1}$ shares one common vertex with $\bigcup_{k=0}^{i}\bbb E_k$, $\bigcup_{k=0}^{i+1}\bbb E_k$ is well-configured. 
If $\bbb E_{i+1}$ is a directed path, its two end-vertices belong to $\bigcup_{k=0}^{i}\bbb E_k$. Since the well-configuration of $\bigcup_{k=0}^{i}\bbb E_k$ guarantees that the two end-vertices have the same value, from Corollary \ref{coro:path}, $\bbb E_{i+1}$ is well-configured, and so is $\bigcup_{k=0}^{i+1}\bbb E_k$. 
By induction, the claim
is established. Since $\bigcup_{k=0}^{p}\bbb E_k = \bbb N$, the proof is complete.
\hfill$\qed$

\vspace{.05in}

The proof of Theorem \ref{thm:steve} provides a constructive approach that systematically designs $C_{ji}$ matrices for a strongly connected multi-agent system to be well-configured.

For each ear decomposition, say $\mathrm{D}=\{\bbb E_0,\bbb E_1,\ldots,\bbb E_p\}$, let $l(\bbb E_i)$ denote the length of ear $\bbb E_i$, i.e., the number of arcs in $\bbb E_i$. Theorem \ref{thm:steve} immediately implies the following sufficient conditions for well-configuration. 

\vspace{.05in}

\begin{corollary}\label{thm:strongsufficient}
Suppose that $\bbb N$ is strongly connected and let $\mathrm{D}$ be an ear decomposition of $\bbb N$. If 
$$\max_{\bbb E\in\mathrm D} l(\bbb E)\le n,$$
then $\bar{\bbb N}$ can be well-configured with all $\scr K_{ji}\neq 0$, $i\in\mathbf{m}$, $j\in\scr N_i$.
\end{corollary}

\vspace{.05in}

More can be said. 
For a strongly connected graph $\bbb G$, write $\scr D$ for the set of all possible ear decompositions of $\bbb G$. 
Define
$$\chi(\bbb G)=\min_{\mathrm D\in\scr D}\max_{\bbb E\in\mathrm D} l(\bbb E).$$
Since each ear decomposition begins with a directed cycle and the shortest possible length of a cycle is two, e.g., a pair of agents which are neighbors of each other, $\chi(\bbb G)\ge 2$.


\vspace{.05in}

\begin{corollary}\label{coro:allear}
If $\bbb N$ is strongly connected and $\chi(\bbb N)\le n$, then $\bar{\bbb N}$ can be well-configured with all $\scr K_{ij}\neq 0$, $i\in\mathbf{m}$, $j\in\scr N_i$. 
\end{corollary}

\vspace{.05in}

Although Corollary \ref{coro:allear} provides a weaker condition, to our knowledge, it is still an open problem to construct an efficient algorithm to find all ear decompositions of a strongly connected graph.


\subsection{Symmetric Directed Graphs}

A directed graph is called {\em symmetric} if whenever $(i,j)$ is an arc in the graph, so is $(j,i)$. 
A symmetric directed graph is often called undirected in the literature, which simplifies each pair of directed edges, say $(i,j)$ and $(j,i)$, to one undirected edge between vertices $i$ and $j$. We stick to the term ``symmetric directed graphs'' because of definition of the incidence matrix given in Section \ref{sec:well system}. Consider a symmetric directed graph with $m$ vertices and $d$ directed edges. Then $d$ must be an even number. Our definition of an incidence matrix is of size $m\times d$, while the standard definition of an incidence matrix of the corresponding undirected graph is of size $m\times (d/2)$. Thus using the term ``undirected'' may cause confusion. 
It is worth noting that rooted and strong connectedness boil down to the same connectivity for symmetric directed graphs. 


For any symmetric directed graph $\bbb G$, since each pair of arcs between any pair of neighboring agents in a symmetric directed graph is a cycle with length 2, all these cycles form an ear decomposition, which leads to $\chi(\bbb G)=2$. The following necessary and sufficient condition on well-configuration for symmetric directed graphs is easy to derive from Corollary \ref{coro:allear}.

\vspace{.05in}

\begin{theorem}\label{thm:symmetric}
If $\bbb N$ is a symmetric directed graph, then $\bar{\bbb N}$ can be well-configured with all $\scr K_{ji}\neq 0$, $i\in\mathbf{m}$, $j\in\scr N_i$, if and only if $\bbb N$ is strongly connected and $n\ge 2$.
\end{theorem}

\vspace{.05in}

As will be seen in the next section, there is a motivation, for the purpose of algorithm design, to figure out a condition under which a symmetric directed graph can be well-configured with the additional constraint that $C_{ij}=C_{ji}$ for all $i\in\mathbf{m}$ and $j\in\scr N_i$. To this end, we need the following modified concept of ear decompositions.

A {\em symmetric ear decomposition} of a symmetric directed graph without self-arcs $\bbb G = (\scr V, \scr E)$ with at least two vertices is a sequence of symmetric subgraphs of $\bbb G$, denoted $\{\bbb E_0,\bbb E_1,\ldots,\bbb E_p\}$, in which $\bbb E_0$ is a symmetric directed cycle, and each $\bbb E_i$, $i\in\mathbf{p}$, is a symmetric directed path or a symmetric directed cycle with the following properties:
\begin{enumerate}
   \item $\{\bbb E_0,\bbb E_1,\ldots,\bbb E_p\}$ form an arc partition of $\bbb G$, i.e., 
   $\bbb E_i$ and $\bbb E_j$ are arc disjoint if $i\neq j$, and $\bigcup_{k=0}^{p}\bbb E_k = \bbb G$;
   \item For each $i\in\mathbf{p}$, if $\bbb E_i$ is a symmetric directed cycle, then it has precisely one vertex in common with $\bigcup_{k=0}^{i-1}\bbb E_k$; if $\bbb E_i$ is a symmetric directed path, then its two end-vertices are the only two vertices in common with $\bigcup_{k=0}^{i-1}\bbb E_k$.
\end{enumerate}
Each of $\bbb E_0,\bbb E_1,\ldots,\bbb E_p$ is called a {\em symmetric ear} of the decomposition.
Not all symmetric directed graphs admit a symmetric ear decomposition. A symmetric directed graph is called {\em $k$-connected} if, upon removal of any $k-1$ two-length cycles, the resulting graph is still strongly connected. It has been proved that a symmetric directed graph has a symmetric ear decomposition if and only if it is 2-connected \cite{undirectedear}.\footnote{This is because a symmetric ear decomposition of a symmetric directed graph is essentially equivalent to an ear decomposition of an undirected graph, and a $k$-connected symmetric directed graph is essentially equivalent to a $k$-edge-connected undirected graph.} 
A 2-connected symmetric directed graph may admit multiple symmetric ear decompositions, and apparently, the number of all possible different symmetric ear decompositions is finite.
For each symmetric ear decomposition, say $\mathrm{D}=\{\bbb E_0,\bbb E_1,\ldots,\bbb E_p\}$, let $l(\bbb E_i)$ denote the length of symmetric ear $\bbb E_i$, i.e., the number of two-length cycles in $\bbb E_i$. 
Using the same arguments as in the proof of Theorem \ref{thm:steve}, we have the following result.

\vspace{.05in}

\begin{theorem}
Suppose that $\bbb N$ is 2-connected symmetric directed graph and let $\mathrm{D}$ be a symmetric ear decomposition of $\bbb N$. If for each symmetric ear $\bbb E\in\mathrm{D}$, $\{\scr K_{ij}f\lor \scr K_{ji}: (i,j)\in\bbb E\}$\footnote{We use $\{a\lor b\}$ to denote that either $a$ or $b$ is an element in the set.} is an independent family, then $\bar{\bbb N}$ is well-configured by matrices $C_{ij}=C_{ji}$, $i\in\mathbf{m}$, $j\in\scr N_i$. If, in addition,  
$\max_{\bbb E\in\mathrm D} l(\bbb E)\le n$,
then $\bar{\bbb N}$ can be well-configured with all $\scr K_{ij}=\scr K_{ji}\neq 0$, $i\in\mathbf{m}$, $j\in\scr N_i$.
\end{theorem}

\vspace{.05in}

In the sequel, we will propose and analyze a few distributed algorithms for well-configured systems under different scenarios. 

\vspace{.1in}

\section{Algorithms for Symmetric Directed Graphs}

In this section, we assume that the neighbor graph is symmetric and $C_{ij}=C_{ji}$ whenever agents $i$ and $j$ are a pair of neighbors. 
We begin with the simplest case in which the neighbor graph is fixed.


\subsection{Fixed Symmetric Directed Graphs}

Consider any strongly connected, symmetric directed graph $\bbb N$ with $m$ agents. Our first algorithm appeals to the idea of gradient descent in convex optimization, which is for each agent $i$,
\begin{align}
    x_i(t+1) &= x_i(t)-\alpha(t) \sum_{j\in\mathcal{N}_i} \Big[(C_{ij}'C_{ij}+ C_{ji}'C_{ji})\nonumber\\ 
    & \;\;\;\;\;\;\;\;\;\;\;\;\;\;\;\;\;\;\;\;\;\;\;\;\;\;\;\;\;\;\;\;\;\;\;\times (x_i(t) - x_j(t))\Big],\label{eq:undirected}
\end{align}
where $\alpha(t)$ is a positive time-varying stepsize satisfying $\sum_{t}\alpha(t)=\infty$ and $\sum_{t}\alpha^2(t)<\infty$.

\vspace{.05in}

\begin{theorem}\label{thm:gradient}
If $\bbb N$ is a strongly connected symmetric directed graph and $\bar{\bbb N}$ is well-configured, then algorithm \eqref{eq:undirected} will lead all the agents to reach a consensus. 
\end{theorem}

\vspace{.05in}

The algorithm \eqref{eq:undirected} involves a term $(C_{ij}'C_{ij}+ C_{ji}'C_{ji})(x_i(t) - x_j(t))$ in each agent $i$'s update, where $j$ is any neighbor of agent $i$. In the case when $C_{ij}\neq C_{ji}$, it will require that each agent $i$ receives two signals, $C_{ji}x_j(t)$ and $C_{ij}x_j(t)$, from each of its neighbors at each time step. Although allowing  $C_{ij}\neq C_{ji}$ in a symmetric directed graph makes well-configuration easier in light of Theorem \ref{thm:symmetric}, transmitting two signals could be an issue in communication. In the case when $C_{ij} = C_{ji}$ so that only one signal is transferred, the underlying symmetric directed graph will need to be 2-connected to guarantee well-configuration. These facts are true for all the remaining algorithms in this section.

\vspace{.05in}

{\bf Proof of Theorem \ref{thm:gradient}:}
The $m$ update equations in \eqref{eq:undirected} can be combined into one state form as
\begin{align*}
    x(t+1) = x(t)-\alpha(t) \bar J C' C\bar J' x(t),
\end{align*}
where $x={\rm column}\{x_1, x_2, \ldots, x_m\}$,
which is exactly the gradient descent of minimizing the convex function $\|C\bar J' x\|_2^2$. Thus with appropriate time-varying stepsize $\alpha(t)$ (i.e., $\sum_{t}\alpha(t)=\infty$ and $\sum_{t}\alpha^2(t)<\infty$), $x(t)$ will asymptotically converge to an optimal point of $\|C\bar J' x\|_2^2$, which must be a consensus vector as ${\rm kernel}\;C\bar J'={\rm span}\;\bar I$.
\hfill$\qed$

\vspace{.05in}

The above algorithm requires all $m$ agents share the same sequence of diminishing stepsizes. Our second algorithm gets around this limitation and is thus fully distributed, which is described as follows.


Since well-configuration only depends on $\scr K_{ij}$, the kernel of $C_{ij}$, $i\in\mathbf{m}$, $j\in\scr N_i$, 
without loss of generality, we assume each $C_{ij}$ has full row rank and its rows are orthonormal, which implies that $C_{ij}C'_{ij}=I$ and $P_{ij}\dfb C'_{ij}(C_{ij}C'_{ij})^{-1}C_{ij}=C'_{ij}C_{ij}$ is an orthogonal projection matrix. For each agent $i\in\mathbf{m}$,
\begin{align}
    x_i(t+1) &= x_i(t)-\frac{1}{2(d_i+1)} \sum_{j\in\mathcal{N}_i} \Big[(C_{ij}'C_{ij}+ C_{ji}'C_{ji})\nonumber\\ 
    & \;\;\;\;\;\;\;\;\;\;\;\;\;\;\;\;\;\;\;\;\;\;\;\;\;\;\;\;\;\;\;\;\;\;\;\times (x_i(t) - x_j(t))\Big].\label{eq:undirected_distributed}
\end{align}


\begin{theorem}\label{thm:fixedstep}
If $\bbb N$ is symmetric, strongly connected and $\bar{\bbb N}$ is well-configured, then algorithm \eqref{eq:undirected_distributed} will lead all the agents to reach a consensus exponentially fast. 
\end{theorem}

\vspace{.05in}

To prove the theorem, we need the following lemmas. 

\vspace{.05in}

\begin{lemma} \label{lem:jccj}
If $\bar{\bbb N}$ is well-configured, then $\bar J C'C\bar J'$ is positive semidefinite with exactly $m$ eigenvalues at zero. 
\end{lemma} 

\vspace{.05in}

{\bf Proof of Lemma \ref{lem:jccj}:}
It is clear that $\bar J C'C\bar J'$ is positive semidefinite. 
Since ${\rm kernel}\;C\bar J'={\rm span}\;\bar I$, 
$\bar JC'C\bar J'$ has exactly rank$(\bar I) = n$ eigenvalues at zero. 
\hfill$\qed$

\vspace{.05in}

\begin{lemma} \label{lem:Wjccj}
Let $\bar W = W\otimes I$, where $W$ is a positive diagonal matrix. If $\bar{\bbb N}$ is well-configured, then $\bar W\bar J C'C\bar J'$ has exactly $n$ eigenvalues at zero, and all the remaining eigenvalues are positive. 
\end{lemma}

\vspace{.05in}

{\bf Proof of Lemma \ref{lem:Wjccj}:}
Note that $\bar W\bar JC'C\bar {J'}$ has the same spectrum as $\bar W^{-\frac{1}{2}}\bar W\bar J C'C\bar J'\bar W^{\frac{1}{2}}=\bar W^{\frac{1}{2}}\bar JC'C\bar J' \bar W^{\frac{1}{2}}$. It is clear that $\bar W^{\frac{1}{2}}\bar JC'C\bar J'\bar W^{\frac{1}{2}}$ is positive semidefinite. 
In addition, from Lemma~\ref{lem:jccj}, $C\bar J' \bar W^{\frac{1}{2}} x = 0$ if and only if $\bar W^{\frac{1}{2}} x \in \text{span}\; \bar I$, which implies that $\text{kernel}\; C\bar J' \bar W^{\frac{1}{2}} = \text{span}\; \bar W^{-\frac{1}{2}} \bar I$ and $\bar W\bar J C'C\bar J'$ contains exactly $n$ eigenvalues at zero.
\hfill$\qed$

\vspace{.05in}

We also need the following ``mixed matrix norm'' concept introduced in \cite{tacle}.
Let    $\|\cdot \|_{\infty}$ denote   the
  induced infinity  norm   and write
 $\R^{mn\times mn}$ for  the vector space of all $m\times m$  block matrices $Q = [Q_{ij}]$
whose $ij$th entry     is an  $n\times n$ matrix $Q_{ij}\in\R^{n\times n}$.
Define the $(2,\infty)$ norm of $Q\in\R^{mn\times mn}$, written $\|Q\|_{2,\infty}$, to be
$$\|Q\|_{2,\infty} = \|\langle Q\rangle \|_{\infty},$$
 where $\langle Q\rangle $ is the $m\times m$ matrix in $\R^{m\times m}$  whose $ij$th entry is
$\|Q_{ij}\|_2$,
 where $\|\cdot\|_2$ denotes the induced 2-norm.
It has been shown in \cite[Lemma 3]{tacle} that $\|\cdot \|_{2,\infty}$ is a sub-multiplicative matrix norm.

In the sequel, for a matrix $Q\in\R^{mn\times mn}$, we sometimes use $[Q]_{ij}$, $i,j\in\mathbf{m}$, to denote the $ij$th block of $Q$, which is an $n\times n$ matrix.

We are now in a position to prove Theorem \ref{thm:fixedstep}.

\vspace{.05in}

{\bf Proof of Theorem \ref{thm:fixedstep}:}
The $m$ update equations in \eqref{eq:undirected_distributed} can be combined into one state form as
\begin{align*}
    x(t+1) = x(t)-\bar D \bar JC' C\bar J' x(t),
\end{align*}
where $\bar D={\rm blockdiag}\{\frac{1}{2(d_1+1)}I_m, \ldots, \frac{1}{2(d_m+1)}I_m\}$.
Consider each block of $\bar D \bar {J}C' C\bar J'$. For any $i\in\mathbf{m}$,  
\begin{align*}
    [\bar D \bar {J}C' C\bar J']_{ii} 
    &= \frac{1}{2(d_i+1)} \sum_{j\in\mathcal{N}_i} (C_{ij}'C_{ij}+C_{ji}'C_{ji}) \\
    &= \frac{1}{2(d_i+1)} \sum_{j\in\mathcal{N}_i} (P_{ij}+P_{ji}),
\end{align*}
and for any $i\in\mathbf{m}$, $j\in\mathcal{N}_i$, 
\begin{align*}
    [\bar D \bar {J}C' C\bar J']_{ij} 
    &= -\frac{1}{2(d_i+1)}(C_{ij}'C_{ij}+C_{ji}'C_{ji}) \\
    &= -\frac{1}{2(d_i+1)} (P_{ij}+P_{ji}).
\end{align*}
Since each $P_{ij}$ is an orthogonal projection matrix, $\| P_{ij} \|_2 = 1$. 
Consider the $2$-norm for each block. For any $i\in\mathbf{m}$,
\begin{align*}
    \|[\bar D \bar {J}C' C\bar J']_{ii}\|_2 
    &\le \frac{1}{2(d_i+1)} \sum_{j\in\mathcal{N}_i} (\| P_{ij}\|_2+\|P_{ji}\|_2) \\
    &\le \frac{d_i}{d_i+1}, 
\end{align*}
and for any $i\in\mathbf{m}$, $j\in\mathcal{N}_i$,
$$\|[\bar D \bar {J}C' C\bar J']_{ij}\|_2 \le \frac{\|P_{ij}\|_2+\|P_{ji}\|_2}{2(d_i+1)}\le \frac{1}{d_i+1}.$$
Next consider the $(2,\infty)$-norm of $\bar D \bar {J}C' C\bar J'$:
\begin{align*}
&\;\;\;\;\;\; \|\bar D \bar {J}C' C\bar J'\|_{2,\infty} \\
&= \min_{i\in\mathbf{m}} \|[\bar D \bar {J}C' C\bar J']_{ii}\|_2 + \sum_{j\in\mathcal{N}_i}\|[\bar D \bar {J}C' C\bar J']_{ij}\|_2\\
&\le \min_{i\in\mathbf{m}} \Big(\frac{d_i}{d_i+1} + \sum_{j\in\mathcal{N}_i}\frac{1}{d_i+1}\Big)
\le \min_{i\in\mathbf{m}} \frac{2d_i}{d_i+1} < 2,
\end{align*} 
which implies that the spectral radius of $\bar D \bar {J}C' C\bar J'$ is less than 2. 
It follows that $I-\bar D \bar {J}C' C\bar J'$ has $n$ eigenvalues at one and all the other eigenvalues lie in $(-1,1)$, which implies that $x(t)$ will reach a consensus exponentially fast. 
\hfill$\qed$

\subsection{Time-varying Symmetric Directed Graphs}

In this subsection, we consider the following scenario of time-varying symmetric directed graphs. Let an $m$-vertex symmetric directed graph $\bbb N$ represent all allowable communication among the $m$ agents. In other words, agents $i$ and $j$ are allowed to communicate with each other if and only if $(i,j)$ is an arc in $\bbb N$. For each time $t$, we use a time-dependent $m$-vertex symmetric directed graph $\bbb N(t)$ to describe the neighbor relations among the $m$ agents at time $t$. That is, if agents $i$ and $j$ communicate at time $t$, then $(i,j)$ is an arc in $\bbb N(t)$. It is easy to see that $\bbb N(t)$ is a spanning subgraph of $\bbb N$, and all such possible spanning subgraphs is a finite set. We assume that $\bar{\bbb N}$ is well-configured, i.e., each arc $(i,j)$ in $\bbb N$ is associated with a matrix $C_{ij}$  such that ${\rm kernel}\; C\bar J'={\rm span}\; \bar I$, with $J$ being the incidence matrix of~$\bbb N$. 

For any time-varying symmetric directed graph sequence just described, we propose the following algorithm using the Metropolis weights:
\begin{align}
    x_i(t+1) &= x_i(t)-\frac{1}{2} \sum_{j\in\mathcal{N}_i(t)} \Big[w_{ij}(t)(C_{ij}'C_{ij}+ C_{ji}'C_{ji})\nonumber\\ 
    & \;\;\;\;\;\;\;\;\;\;\;\;\;\;\;\;\;\;\;\;\;\;\;\;\;\;\;\;\;\;\;\;\;\;\;\times (x_i(t) - x_j(t))\Big],\label{eq:metrotimevary}
\end{align}
where $\scr N_i(t)$ is the neighbor set of agent $i$ at time $t$ 
and $w_{ij}(t)$ are the Metropolis weights corresponding to $\bbb N(t)$, which are proposed in \cite{metro2} for solving the distributed averaging problem over symmetric directed graphs and defined as 
$$w_{ij}(t) = \frac{1}{1+\max\{d_i(t), d_j(t)\}}, \;\;\;\;\; j\in\mathcal{N}_{i}(t),$$
where $d_i(t) = |\mathcal{N}_i(t)|$ denotes the number of neighbors of agent $i$ at time $t$. 

\vspace{.05in}

\begin{theorem}\label{thm:varyundirected}
Suppose that $\bar{\bbb N}$ is well-configured. 
If $\bbb N$ is symmetric, strongly connected and each edge of $\bbb N$ appears infinitely often in the infinite sequence of neighbor graphs $\bbb N(1),\bbb N(2),\bbb N(3),\ldots$, then algorithm \eqref{eq:metrotimevary} will guarantee all $m$ agents to reach a consensus. 
\end{theorem}

\vspace{.05in}

To prove the theorem, we first combine the $m$ update equations in \eqref{eq:metrotimevary} into one state form.
To this end, we tailor the definition of an incidence matrix for spanning subgraphs as follows. Consider a directed graph $\bbb G$ with $m$ vertices and $d$ directed edges. Let $\scr E$ denote the arc set of $\bbb G$ and $J$ denote the $m\times d$ incidence matrix of $\bbb G$ according to some ordering of the arcs in $\scr E$. Let
$\bbb H$ be a spanning subgraph of $\bbb G$. We define the {\em spanning incidence matrix} of $\bbb H$ as an $m\times d$ matrix in which column $k$ has exactly one $1$ in row $i$ and exactly one $-1$ is row $j$ if the $k$th arc in $\bbb G$ is $(j, i)$ and $(j, i)$ is also an arc in $\bbb H$. It is clear that the spanning incidence matrix of any spanning subgraph of $\bbb G$ has the same size as the incidence matrix of $\bbb G$. If the $k$th arc in $\bbb G$ is not in a spanning subgraph, then the $k$th column of the incidence matrix of $\bbb G$ is replaced by a zero vector in the spanning incidence matrix. 

We also need the following definition. Consider a symmetric directed graph $\bbb G$ with $d$ arcs. Let $\bar{\bbb G}$ be a spanning subgraph of $\bbb G$ which is also symmetric. Since $\bar{\bbb G}$ is symmetric, its Metropolis weights $\bar w_{ij}$ are well-defined; specifically, $\bar w_{ij}=1/(1+\max\{\bar d_i,\bar d_j\})$, where $\bar d_k$ denotes the number of neighbors of vertex $k$ in $\bar{\bbb G}$. Given an ordering of all the arcs in $\bbb G$, the {\em spanning weight matrix} of $\bar{\bbb G}$ is the $d\times d$ diagonal matrix whose $k$th diagonal entry equals $\bar w_{ij}$ if the $k$th arc in $\bbb G$ is $(j,i)$ and $(j,i)$ is also an arc in $\bar{\bbb G}$, or 0 if the $k$th arc in $\bbb G$ is not in $\bar{\bbb G}$.

With the above definitions, it is not hard to verify that the $m$ update equations in \eqref{eq:metrotimevary} can be written as 
\begin{align}\label{eq:metrotimevarysystem}
    x(t+1) = x(t)-\frac{1}{2}\bar J(t)C'\bar W(t) C\bar J'(t) x(t),
\end{align}
where $\bar J(t)=J(t)\otimes I_n$ with  $J(t)$ being the spanning incidence matrix of $\bbb N(t)$, and $\bar W(t)=W(t)\otimes I_n$ with $W(t)$ being the spanning weight matrix of $\bbb N(t)$. It is worth noting that all $W(t)$ are nonnegative diagonal matrices with the same size. It is also worth emphasizing that the definitions of $C$, $J(t)$, and $W(t)$ are based on the same ordering of the arcs in $\bbb N$, and the equality \eqref{eq:metrotimevarysystem} is independent of the ordering. 


To proceed, we need the following concept and result. 

A squre matrix $M$ is called {\em paracontracting} with respect to a vector norm $\|\cdot\|$ if
$\|Mx\|\leq \|x\|$ and the strict inequality holds whenever $Mx\neq x$. 

It is easy to see that any symmetric matrix is paracontracting with respect to the $2$-norm if all its eigenvalues lie in the interval $(-1,1]$. 

For a square matrix $M$, we define its fixed point set as
$$\scr{F}(M) = \left\{x :  Mx=x \right\}.$$
Paracontracting matrices have the following properties. 



\vspace{.05in}

\begin{lemma}
Suppose that a finite set of square matrices $\{M_1,M_2,\ldots,M_p\}$ are paracontracting with respect
to the same vector norm. Let $\sigma(1),\sigma(2),\ldots$ be an infinite sequence of integers
taking values in $\{1,2,\ldots,p\}$ and $\scr{I}$ be the set of all integers that
appears infinitely often in the sequence. Then for any initial vector $z(0)$,
the sequence of vectors generated by $z(t+1)=M_{\sigma(t)}z(t)$ has a limit $z^*\in\bigcap_{i\in\scr{I}}\scr{F}(M_i)$.
\label{key}\end{lemma}

\vspace{.05in}

The lemma is a special case of Theorem~1 in \cite{paracontract}.

We also need the following lemmas. 

\vspace{.05in}

\begin{lemma} \label{lem:jcwcj}
Let $\bar W = W\otimes I$, where $W$ is a positive diagonal matrix. If $\bar{\bbb N}$ is well-configured, then $\bar JC'\bar WC\bar J'$ has exactly $n$ eigenvalues at zero, and all the remaining eigenvalues are positive. 
\end{lemma}

\vspace{.05in}


{\bf Proof of Lemma \ref{lem:jcwcj}:}
It is clear that $\bar JC'\bar WC\bar J'$ is positive semidefinite. 
In addition, for any $x$, there exists a $y$ such that $\bar W^{1/2}C\bar J' x = y$, i.e., $C\bar J' x = W^{-1/2}y$.
From Lemma~\ref{lem:jccj}, we know if and only if $x$ satisfies $x \in \text{span}\; \bar I$, we have $C\bar J' x = 0 = W^{-1/2}y$, which implies that $W^{1/2}C\bar J' x = y=0$, i.e., $\text{kernel}\; W^{1/2} C\bar J' = \text{span}\; \bar I$, and $\bar JC'\bar WC\bar J'$ contains exactly $n$ eigenvalues at zero.
\hfill$\qed$

\vspace{.05in}

\begin{lemma}\label{lm:updatematrix}
Let $\bbb G$ be a symmetric, spanning subgraph of $\bbb N$, $W$ be the spanning weight matrix of $\bbb G$, and $J$ be the spanning incidence matrix of $\bbb G$. Then all the eigenvalues of $I-\frac{1}{2}\bar J C'\bar W C\bar J'$ lie in $(-1,1]$. If furthermore $\bbb G=\bbb N$, $I-\frac{1}{2}\bar J C'\bar W C\bar J'$ has exactly $n$ eigenvalues at one and all the remaining eigenvalues lie in $(-1,1)$.
\end{lemma}

\vspace{.05in}

{\bf Proof of Lemma \ref{lm:updatematrix}:}
Consider each block of $\frac{1}{2} \bar J C'\bar W C\bar J'$. 
For any $i\in\mathbf{m}$, 
\begin{align*}
    [ \frac{1}{2}\bar J C'\bar W C\bar J']_{ii} 
    &= \frac{1}{2} \sum_{j\in\mathcal{N}_i} (w_{ij} C_{ij}'C_{ij}+w_{ji}C_{ji}'C_{ji}) \\
    &= \frac{1}{2} \sum_{j\in\mathcal{N}_i} (w_{ij}P_{ij}+w_{ji}P_{ji}),
\end{align*}
and for any $i\in\mathbf{m}$, $j\in\mathcal{N}_i$,
\begin{align*}
    [ \frac{1}{2}\bar J C'\bar W C\bar J']_{ij} &= -\frac{w_{ij}}{2}(C_{ij}'C_{ij}+C_{ji}'C_{ji}) \\
    &= -\frac{w_{ij}}{2} (P_{ij}+P_{ji}).
\end{align*}
Note that $\| P_{ij} \|_2 = 1$, $w_{ij}=w_{ji}\le 1$, and $\sum_{j\in\mathcal{N}_i} w_{ij} < 1$.
Considering the $2$-norm for each block. For any $i\in\mathbf{m}$,
\begin{align*}
    \|[ \frac{1}{2} \bar J C'\bar W C\bar J']_{ii}\|_2 &\le \frac{1}{2} \sum_{j\in\mathcal{N}_i} w_{ij} ( \| P_{ij}\|_2+\|P_{ji}\|_2)\\ &\le \sum_{j\in\mathcal{N}_i}w_{ij}, 
\end{align*}
and for any $i\in\mathbf{m}$, $j\in\mathcal{N}_i$,
$$\|[\frac{1}{2}\bar J C'\bar W C\bar J']_{ij}\|_2 \le \frac{1}{2} w_{ij}(\|P_{ij}\|_2+\|P_{ji}\|_2)\le w_{ij}.$$
Next consider the $(2,\infty)$-norm of $ \frac{1}{2} \bar J C'\bar W C\bar J'$:
\begin{align*}
    &\;\;\;\;\; \|\frac{1}{2}\bar J C'\bar W C\bar J'\|_{2,\infty} \\
    &= \min_{i\in\mathbf{m}} \|[ \frac{1}{2}\bar J C'\bar W C\bar J']_{ii}\|_2 + \sum_{j\in\mathcal{N}_i}\|[ \frac{1}{2}\bar J C'\bar W C\bar J']_{ij}\|_2 \\
    &\le \min_{i\in\mathbf{m}} \Big(\sum_{j\in\mathcal{N}_i}w_{ij} + \sum_{j\in\mathcal{N}_i}w_{ij}\Big) \\
    &= 2\min_{i\in\mathbf{m}} \sum_{j\in\mathcal{N}_i} \frac{1}{1+\max(d_i, d_j)}
    \le 2\min_{i\in\mathbf{m}}\frac{d_i}{1+\min_{k\in\mathbf{m}}d_k}\\
    &= \frac{2\min_{i\in\mathbf{m}}d_i}{1+\min_{k\in\mathbf{m}}d_k}< 2,
\end{align*} 
which implies that the spectral radius of $\frac{1}{2}\bar J C'\bar W C\bar J'$ is less than 2. 
It follows that $I-\frac{1}{2}\bar J C'\bar W C\bar J'$ has $n$ eigenvalues at one and all the other eigenvalues lie in $(-1,1)$. 


In the case when $\bbb G=\bbb N$, $W$ is positive diagonal matrix. Then the lemma is true by Lemma \ref{lem:jcwcj}. \hfill$\qed$

\vspace{.05in}

The above lemma implies that each update matrix $(I-\frac{1}{2}\bar J(t)C'\bar W(t) C\bar J'(t))$ in \eqref{eq:metrotimevarysystem} is paracontracting with respect to the $2$-norm. 


\vspace{.05in}





\begin{lemma} \label{lem:cj_cWj}
Let $\bbb G_1,\bbb G_2, \ldots,\bbb G_p$ be a finite set of symmetric, spanning subgraphs of $\bbb G$. If the union of  $\bbb G_1,\bbb G_2, \ldots,\bbb G_p$ is $\bbb G$, then ${\rm kernel}\; C\bar J'={\rm kernel}\; C(\sum_{i=1}^p \bar W_i^{1/2}\bar J'_i)$, 
where $J$ is the incidence matrix of $\bbb G$, $J_i$ is the spanning incidence matrix of $\bbb G_i$, and $W_i$ is the spanning weight matrix of $\bbb G_i$.
\end{lemma}

\vspace{.05in}

{\bf Proof of Lemma \ref{lem:cj_cWj}:}
If $(i,j)$ is an edge in $\bbb G$ but not in $\bbb G_k$, then the corresponding Metropolis weight $w_{ij} = 0$ for $\bbb G_k$, which implies that $\bar W_k^{1/2}\bar J'_k = \bar W_k^{1/2}\bar J'$. 
Then $C(\sum_{i=1}^p \bar W_i^{1/2}\bar J'_i) = C(\sum_{i=1}^p \bar W_i^{1/2}) \bar J' = (\sum_{i=1}^p \bar W_i^{1/2}) C \bar J'$. Since the union of  $\bbb G_1,\bbb G_2, \ldots,\bbb G_p$ is $\bbb G$,  $\sum_{i=1}^p \bar W_i^{1/2}$ is a positive diagonal matrix and thus nonsingular, then ${\rm kernel}\; C(\sum_{i=1}^p \bar W_i^{1/2}\bar J'_i) = {\rm kernel}\; (\sum_{i=1}^p \bar W_i^{1/2})C\bar J'= {\rm kernel}\; C\bar J'.$
\hfill$\qed$


\vspace{.05in}

Now we are in a position to prove Theorem \ref{thm:varyundirected}.

\vspace{.05in}

{\bf Proof of Theorem \ref{thm:varyundirected}:}
Let $\scr S$ denote the set of all possible spanning subgraphs of $\bbb N$, which apparently is a finite set. Let $\scr I\subset \scr S$ denote the set of those spanning subgraphs which appears infinitely often in the infinite sequence $\bbb N(1),\bbb N(2),\bbb N(3),\ldots$. 
Denote all spanning graphs in $\scr I$ as $\bbb N_1,\bbb N_2,\ldots,\bbb N_p$. 
Since each edge of $\bbb N$ appears infinitely often, the union of $\bbb N_1,\bbb N_2,\ldots,\bbb N_p$ is $\bbb N$. 

Let $J$ be the incidence matrix of $\bbb N$ and $J_i$ be the spanning incidence matrix of $\bbb N_i$ for all $i\in\mathbf{p}$. Let $W$ and $W_i$ be the spanning weight matrices of $\bbb N$ and $\bbb N_i$, $i\in\mathbf{p}$, respectively.
From Lemma \ref{lm:updatematrix}, each update matrix $(I-\frac{1}{2}\bar J(t)C'\bar W(t) C\bar J'(t))$ in \eqref{eq:metrotimevarysystem} is paracontracting with respect to $2$-norm for all $t$.  
From Lemma \ref{key}, $x(t)$ will asymptotically converge to a common fixed point all $(I-\frac{1}{2}\bar J_iC'\bar W_i C\bar J'_i)$, $i\in\mathbf{p}$.
It is easy to see that $\scr F(I-\frac{1}{2}\bar J_iC'\bar W_i C\bar J'_i)={\rm kernel}\; \bar W_i^{\frac{1}{2}}C\bar J'_i={\rm kernel}\; C\bar W_i^{\frac{1}{2}}\bar J'_i$. Thus $x(t)$ will converge to a point in the intersection of ${\rm kernel}\; C\bar W_i^{\frac{1}{2}}\bar J'_i$, $i\in\mathbf{p}$. 

It is clear that the intersection of ${\rm kernel}\; C\bar W_i^{\frac{1}{2}}\bar J'_i$, $i\in\mathbf{p}$ is a subset of ${\rm kernel}\; C(\sum_{i=1}^p \bar W_i^{1/2}\bar J'_i)$. From Lemma \ref{lem:cj_cWj}, ${\rm kernel}\; C(\sum_{i=1}^p \bar W_i^{1/2}\bar J'_i)={\rm kernel}\; C\bar J'$. Since $\bar{\bbb N}$ is well-configured, ${\rm kernel}\; C\bar J'={\rm span}\; \bar I$, which implies that the intersection of ${\rm kernel}\; \bar CW_i^{\frac{1}{2}}\bar J'_i$, $i\in\mathbf{p}$, is a subset of ${\rm span}\; \bar I$.~\hfill$\qed$






\vspace{.1in}

\section{Algorithms for Directed Graphs}

In this section, we discuss some special types of strongly connected graphs. We begin with directed cycles, the simplest strongly connected graphs.

\subsection{Directed Cycles with Specific Initial States}

Consider an $m$-vertex directed cycle $1\rightarrow 2\rightarrow\cdots\rightarrow m\rightarrow1$ whose local agreement equations are given in \eqref{eq:cycle1}.
The agents update their states as follows:
\begin{align}
    x_i(t+1)&=x_i(t)-\frac{1}{2}P_i(x_i(t)-x_{i-1}(t)),\;\;\; i\in\mathbf{m},\label{eq:directed_update}
\end{align}
where $P_i=C'_i(C_iC'_i)^{-1}C_i$ is a projection on $\scr K_i^\perp$. In this subsection, we assume that each agent $i$ initializes its state $x_i(0)\in {\rm image}\; P_i$, which can be implemented in a distributed manner. 

\vspace{.05in}

\begin{proposition}\label{prop:cycleinitial}
If $\bbb N$ is an $m$-vertex directed cycle, then algorithm \eqref{eq:directed_update} with $x_i(0)\in {\rm image}\; P_i$, $i\in\mathbf{m}$, will lead all $m$ agents to reach a consensus exponentially fast. 
\end{proposition}

\vspace{.05in}

{\bf Proof of Proposition \ref{prop:cycleinitial}:}
It is easy to see that with the initialization  $x_i(t)\in{\rm image}\; P_i$ for all $t$. Then the updates can be written as 
\begin{align*}
    x_i(t+1)&=P_ix_i(t)-\frac{1}{2}P_i(x_i(t)-x_{i-1}(t))\\
    &=P_i\Big(\frac{1}{2}x_i(t)+\frac{1}{2}x_{i-1}(t)\Big),\;\;\;\;\; i\in\mathbf{m},
\end{align*}
which leads to the system update as
$$x(t+1)=P(F\otimes I)x(t),$$
where $P$ is the block diagonal matrix of all $P_i$ and $F$ is the flocking matrix\footnote{The flocking matrix of a directed graph $\bbb G$ is defined as $D_{\bbb G}^{-1}A'_{\bbb G}$, where $D_{\bbb G}$ is the diagonal matrix
whose $i$th diagonal entry is the in-degree of vertex $i$ in $\bbb G$ and $A_{\bbb G}$ is the adjacency matrix of $\bbb G$. A flocking matrix is a stochastic matrix \cite{reachingp1}.} of the cycle.
The update has the same form as the distributed linear equation solver in \cite{tacle}, which guarantees exponentially fast consensus. 
\hfill$\qed$



\subsection{Directed Cycles with Arbitrary Initial States}

In this subsection, we consider directed cycles and algorithm \eqref{eq:directed_update} without any specific initialization.

\vspace{.05in}

\begin{theorem}\label{thm:directedcyclewitharbitraryintialvalue}
If $\bbb N$ is an $m$-vertex directed cycle and $\bar{\bbb N}$ is well-configured, then algorithm \eqref{eq:directed_update} will lead all $m$ agents to reach a consensus exponentially fast for any initial states. 
\end{theorem}

\vspace{.05in}

To prove the theorem, we first rewrite the $m$ equations in \eqref{eq:directed_update} as one state form as
\eq{x(t+1)=Mx(t),\label{eq:cyclesystem}} 
where $M$ is an $mn\times mn$ matrix whose blocks can be easily figured out via \eqref{eq:directed_update}. 
The system update matrix $M$ has the following properties.

\vspace{.05in}

\begin{lemma}\label{lem:lambda<1}
If $\lambda\neq 1$ is an eigenvalue of $M$, then $|\lambda| < 1.$
\end{lemma}

\vspace{.05in}

{\bf Proof of Lemma \ref{lem:lambda<1}:}
Let $v=
[v_1\;\cdots\;v_m]'\neq 0$ be an eigenvector of $M$ for eigenvalue $\lambda$. Then
\begin{align*}
    \Big(I-\frac{P_i}{2}\Big)v_i + \frac{P_i}{2}v_l = \lambda v_i,
\end{align*}
where $l=i+1$ when $i\in\{1,\ldots,m-1\}$, and $l=1$ when $i = m.$
Re-arranging the equation, we have 
\begin{align}\label{eigeneq}
    \Big(\frac{1}{2}P_i-(1-\lambda)I\Big)v_i = \frac{P_i}{2}v_l,
\end{align}
For all $i,j\in\mathbf{m},$ $v_i$ can be divided into $v_i = \alpha_{i,j} + \beta_{i,j},$ where $\alpha_{i,j}\in{\rm image}\;P_j$ and $\beta_{i,j}\in{\rm kernel}\; P_j.$
Since $P_j$ is symmetric, such decomposition is unique for any $j\in\mathbf{m}$ because ${\rm image}\;P_j\oplus{\rm kernel}\; P_j = \R^n$. Thus $\alpha_{i,j}\perp\beta_{h,j}$ for any $i,j,h.$ 
Substituting $v_i = \alpha_{i,i} + \beta_{i,i}$ and $v_l=\alpha_{l,i}+\beta_{l,i}$ to \eqref{eigeneq}, we have
\begin{align*}
    (\lambda-\frac{1}{2})\alpha_{i,i} - (1-\lambda)\beta_{i,i} = \frac{1}{2}\alpha_{l,i}.
\end{align*}
Let $\beta^*_{i,i}$ be the conjugate transpose of $\beta_{i,i}$.
Left product $\beta^*_{i,i}$ on both sides, we have $(1-\lambda)\beta^*_{i,i}\beta_{i,i}=0.$ Since $\lambda\neq 1,$ we have $\beta_{i,i} = 0,$ which implies that $v_i\in{\rm image}\;P_i$ when $\lambda\neq 1.$ Applying this result to \eqref{eigeneq}, we obtain that
\begin{align}\label{verifylambda<1}
    \frac{1}{2}(v_i+\alpha_{l,i}) = \lambda v_i.
\end{align}
It can be shown that there exists an $i\in\mathbf{m},$ such that $v_i\neq v_l$ and $\|v_i\|_2\ge \|v_l\|_2$. Otherwise, for each $i\in\mathbf{m}$, $v_i$ satisfies either (case 1) $\|v_i\|_2< \|v_l\|_2$ or (case 2) $v_i = v_l$. 
If there exists at least an $i$ satisfying case 1,  then $\sum_{i\in\mathbf{m}}\|v_i\|_2 < \sum_{l\in\mathbf{m}}\|v_l\|_2 = \sum_{i\in\mathbf{m}}\|v_i\|_2,$ which contradicts. While if case 2 holds for all $i,$  we obtain from \eqref{eigeneq} that either $\lambda=1$ or $v = 0$, which is not consistent with $\lambda\neq 1$ and $v$ being nonzero.

For such $i$ above, combining the results that $\|v_{l}\|^2_2 = \|\alpha_{l,i}\|^2_2 + \|\beta_{l,i}\|^2_2$ and that $\|v_i\|_2\ge \|v_l\|_2$ together, we have $\|\alpha_{l,i}\|_2\le\|v_i\|_2.$ It can be shown that $v_i\neq\alpha_{l,i},$ as otherwise if $v_i=\alpha_{l,i},$ then from the fact that $\|\alpha_{l,i}\|_2\le\|v_i\|_2,$ we obtain that $\beta_{li} = 0,$ further implying that $v_i = \alpha_{l,i} = v_{l},$ which is not consistent with $v_i\neq v_l.$ Since $v_i\neq\alpha_{l,i}$ and $\|\alpha_{l,i}\|_2\le\|v_i\|_2$ together imply that $\langle v_i,\alpha_{l,i}\rangle$ is strictly less than $\|v_i\|^2$. After taking two-norm on both side of the \eqref{verifylambda<1}, we have
\begin{align*}
    |\lambda|^2\|v_i\|^2 &= \frac{1}{4}(\|v_i\|^2+\|\alpha_{l,i}\|^2+2\langle v_i,\alpha_{l,i}\rangle)<\|v_i\|^2,
\end{align*}
which implies $|\lambda|<1.$
\hfill$\qed$

\vspace{.05in}



\begin{lemma}\label{lem:fixedpint=spanI}
If $\bar{\bbb N}$ is well-configured, $\{x : Mx=x\} = {\rm span}\; \bar I$.
\end{lemma}

\vspace{.05in}

{\bf Proof of Lemma \ref{lem:fixedpint=spanI}:}
It is easy to verify that $\{x : Mx=x\} \supset {\rm span}\; \bar I,$ thus we focus on $\{x : Mx=x\} \subset {\rm span}\; \bar I$ in the remaining proof.

Let $x$ be any eigenvector for eigenvalue 1. Then for each $i\in\mathbf{m},$
\begin{align*}
    \Big(I-\frac{P_i}{2}\Big)x_i + \frac{P_i}{2}x_l =  x_i,
\end{align*}
where $l=i+1$ when $i\in\{1,\ldots,m-1\}$, and $l=1$ when $i = m.$
Re-arranging the equation, we obtain that $P_i(x_i-x_l)=0.$ For any $i\in\mathbf{m},$ let $y_i = x_i-x_l,$ we have $y_i\in\ker P_i = \scr K_i.$ Taking summation on both sides, then\[\sum_{i\in\mathbf{m}}y_i = \sum_{i\in\mathbf{m}}x_i - \sum_{l\in\mathbf{m}}x_l = 0.\] 
Since for any  $i\in\mathbf{m}$, $y_i = 0$ whenever $\scr K_i = 0$, we have
\eq{\sum_{i\in\mathbf{m},\;\scr K_i\neq 0}y_i = 0.\label{eq:xxx}}
Since $\bar{\bbb N}$ is well-configured, from Lemma~\ref{lem:cycle}, $\{\scr K_i : i\in\mathbf{m}, \scr K_i \neq 0\}$ is an independent family. 
Since $y_i\in\scr K_i$ for each $i\in\mathbf{m}$, \eqref{eq:xxx} implies that $y_i = 0$ for all $i\in\mathbf{m}$ such that $\scr K_i\neq 0$. Thus $y_i = 0$ for all $i\in\mathbf{m}$, namely, $x_i$ all equal, which implies $x\in {\rm span}\; \bar I,$ thus $\{x : Mx=x\} \subset {\rm span}\; \bar I.$
\hfill$\qed$

\vspace{.05in}

\begin{lemma}\label{lem:lambda=1}
If $\bar{\bbb N}$ is well-configured, then $M$ has exactly $n$ eigenvalues at one. 
\end{lemma}


\vspace{.05in}

{\bf Proof of Lemma \ref{lem:lambda=1}:}
To prove the lemma, it is sufficient to show that the algebraic and geometric multiplicity of eigenvalue one are equal, namely, all the Jordan blocks of eigenvalue one are of size 1. This and Lemma~\ref{lem:fixedpint=spanI} imply that $M$ has exactly $n$ eigenvalues at one.

Consider the Jordan form of $M,$ since $|\ker (I-M)^2| - |\ker (I-M)|$ equals the number of Jordan blocks of eigenvalue 1 with size larger than 1, and $\ker (I-M)^2\supset \ker (I-M),$ then to prove all the Jordan blocks of eigenvalue one are of size 1, i.e., $|\ker (I-M)^2| - |\ker (I-M)|=0$, we only need to prove that $\ker (I-M) = \ker (I-M)^2$.

Suppose by contradiction, there exists a nonzero $y\in\ker (I-M)^2\setminus\ker (I-M),$ then  $(I-M)y\in\ker (I-M)$ while $y\notin \ker (I-M).$ Since from Lemma~\ref{lem:fixedpint=spanI}, $\ker(I-M) = \text{span }\bar I, $ there exists a nonzero vector $z,$ such that $(I-M)y = \1\otimes z.$ Expanding the equation, we have $\frac{P_i}{2}(y_i-y_l) = z,$
where $l = i+1$ when $i\in\{1,\ldots,m-1\}$ and $l = 1$ when $i = m.$ This implies $z\in\bigcap_{i\in\mathbf{m}}{\rm image}\;P_i,$ and that
\[0 = \frac{P_i}{2}(y_i-y_l) - z = \frac{P_i}{2}(y_i-y_l-2z),\] then we obtain that
\begin{align*}
    y_i - y_l -2z \in \ker P_i.
\end{align*}
After taking summation on both sides, we have 
\[-2mz = \sum_{i\in\mathbf{m}}y_i-\sum_{l\in\mathbf{m}}y_l-2mz \in \sum_{i\in\mathbf{m}}\ker P_i,\]
so $z\in\sum_{i\in\mathbf{m}}\ker P_i.$ While since $z\in\bigcap_{i\in\mathbf{m}}{\rm image}\;P_i,$ we have $z\perp\ker P_i$ for all $i\in\mathbf{m}$ as $P_i$ is symmetric, thus $z\perp\sum_{i\in\mathbf{m}}\ker P_i,$ this along with $z\in\sum_{i\in\mathbf{m}}\ker P_i$ implies that $z = 0,$  which is not consistent with the assumption that $z$ is a nonzero vector. This way we complete the proof.
\hfill$\qed$

\vspace{.05in}

We are now in a position to prove Theorem \ref{thm:directedcyclewitharbitraryintialvalue}.

\vspace{.05in}

{\bf Proof of Theorem \ref{thm:directedcyclewitharbitraryintialvalue}:}
From Lemmas~\ref{lem:lambda<1}--\ref{lem:lambda=1}, the linear system \eqref{eq:cyclesystem} will converge to the eigenspace of eigenvalue one as $t\rightarrow\infty$. From Lemma \ref{lem:fixedpint=spanI}, the eigenspace of eigenvalue one is ${\rm span}\; \bar I$, which implies that all $x_i(t)$, $i\in\mathbf{m}$, will reach a consensus. Since the linear system is time-invariant, the consensus will be reached exponentially fast. 
\hfill$\qed$

\subsection{A Counterexample}

One may conjecture that following algorithm 
\eq{x_i(t+1)=x_i(t)-\frac{1}{d_i+1}\sum_{j\in\scr N_i}P_{ji}(x_i(t)-x_j(t)),\label{eq:steve}}
where $P_{ji}=C'_{ji}(C_{ji}C'_{ji})^{-1}C_{ji}$ is a projection on $\scr K_{ji}^\perp$,
will lead to a consensus for any strongly connected graphs, considering algorithm \eqref{eq:directed_update}
is a special case of \eqref{eq:steve}.
It turns out that it is not the case, as shown in the following counterexample.

Consider a strongly connected graph with 3 vertices and 4 directed edges $(1,2)(2,3)(3,1)(2,1)$.
For simplicity, we write 4 matrices as $C_1,C_2,C_3,C_4$ for the 4 directed edges whose kernels are $\scr K_1,\scr K_2,\scr K_3,\scr K_4$. The corresponding local agreement equations are
\begin{align*}
    C_{1}(x_2-x_1)&=0\\
    C_{2}(x_3-x_2)&= 0\\
    C_{3}(x_1-x_3)&=0\\
    C_{4}(x_1-x_2)&=0
\end{align*}
{\bf Claim:}
The weighted neighbor graph $\bar{\bbb N}$ is well-configured with all $\scr K_i\neq 0$ if and only if $\{\scr K_1\cap\scr K_4,\scr K_2,\scr K_3\}$ is an independent family, i.e.,
\begin{align} \label{eq:3agentiff}
    &\scr K_1 \cap \scr K_4 \cap (\scr K_2+\scr K_3) =0 \nonumber\\
    & \scr K_2 \cap (\scr K_1 \cap \scr K_4+\scr K_3) =0\\
    & \scr K_3 \cap (\scr K_1 \cap \scr K_4+\scr K_2)=0.\nonumber
\end{align}


To prove the claim,
let $p_1 = x_1-x_2$, $p_2=x_2-x_3$ and $p_3=x_3-x_1$. From local agreement 
\begin{align*}
    C_{1}(x_1-x_2)&=C_{1}p_1=0\\
    C_{2}(x_2-x_3)&= C_{2}p_2=0\\
    C_{3}(x_3-x_1)&=C_{3}p_3=0\\
    C_{4}(x_1-x_2)&=C_{4}p_1=0\\
    p_1+p_2+p_3&=0,
\end{align*}
from which
\begin{align*}
    p_1 &\in \scr K_1 \cap \scr K_4 \cap (\scr K_2+\scr K_3)\\
    p_2 &\in \scr K_2 \cap (\scr K_1 \cap \scr K_4+\scr K_3)\\
    p_3 &\in \scr K_3 \cap (\scr K_1 \cap \scr K_4+\scr K_2).
\end{align*}

We first prove the sufficiency. Suppose \eqref{eq:3agentiff} holds, that is, $\{ \scr K_1 \cap \scr K_4 , \scr K_2, \scr K_3 \}$ is an independent family, which implies that $p_1 = p_2 = p_3 = 0$, i.e., $x_1 = x_2 = x_3.$

We next prove the necessity. 
Suppose to the contrary that $\bar{\bbb N}$ is well-configured but \eqref{eq:3agentiff} does not hold. Then there must exist a nonzero vector $y_1 \in \scr K_1 \cap \scr K_4 \cap (\scr K_2+\scr K_3)$, which implies that there exist $y_2\in \scr K_2$ and $y_3\in\scr K_3$ such that $y_1 = y_2+y_3$ and $y_1\in \scr K_1 \cap \scr K_4$. Since $y_1$ is nonzero, so is either $y_2$ or $y_3$. Let $x_1 = y_1$, $x_2 = 0 $ and $x_3=y_2$, it is easy to check that all the local agreement equations hold:
\begin{align*}
    C_{1}(x_1-x_2)&=C_{1}y_1=0\\
    C_{2}(x_2-x_3)&= -C_{2}y_2=0\\
    C_{3}(x_3-x_1)&= -C_{3}y_3=0\\
    C_{4}(x_1-x_2)&=C_{4}y_1=0\\
    p_1+p_2+p_3&=0,
\end{align*}
while, since $x_1 \neq x_2$, $\bar{\bbb N}$ is not well configured. 
This completes the proof.
\hfill$\qed$

\vspace{.05in}

For this example, the update matrix 
$$M=\begin{bmatrix}
I-\frac{1}{3}P_4-\frac{1}{3}P_3 & \frac{1}{3}P_4 & \frac{1}{3}P_3 \\
\frac{1}{2}P_1 & I-\frac{1}{2}P_1 & 0 \\
0 & \frac{1}{2}P_2 & I-\frac{1}{2}P_2 
\end{bmatrix}
$$
where $P_i=C'_i(C_iC'_i)^{-1}C_i$. However, its eigenspace of eigenvalue one can be larger than ${\rm span}\; \bar I$ even when $\bar{\bbb N}$ is well-configured. To see this, set $C_1 = C_2$, $C_3 = C_4$ and $\scr K_2 \cap \scr K_3=0$, which implies that $\{\scr K_1\cap\scr K_4,\scr K_2,\scr K_3\}$ is an independent family, and thus $\bar{\bbb N}$ is well-configured. Pick any nonzero $y\in\scr K_1$. Then it is easy to verify that 
\[x = \begin{bmatrix}0\\y\\-y\end{bmatrix}\in\left\{x : Mx=x \right\}\]
which implies that ${\rm span}\; \bar I$ is a proper subset of $\left\{x : Mx=x \right\}$. Thus $x(t+1)=Mx(t)$ may converge to a non-consensus state, which has also been validated by simulations. 

Therefore the following two questions remain open. First, what are the graphical conditions on $\bbb N$ under which algorithm \eqref{eq:steve} will lead all the agents to reach a consensus for arbitrary initial states? Second, how one can construct a distributed algorithm for each agent which will drive the system from arbitrary start
states to a consensus for any strongly connected graphs?

\vspace{.1in}

\section{Conclusion}

In this paper, we have studied the problem of achieving a consensus in the face of limited
information transfer between agents, in which each agent receives a linear function of the state of each of its neighbors; in the case when the linear function is realized by a matrix whose kernel is nonzero, the neighbor's state cannot be determined by the information transferred. From this perspective, the problem studied here is related to so-called privacy preserving consensus problems \cite{private}, which typically rely on carefully designed additive noise. The limited information idea here can be used to protect the privacy of agents' states without adding noise. 
The problem is also related to the compressed communication techniques which have been recently used to address the communication bottleneck in distributed optimization and machine learning \cite{compressed}.

The feasibility of the problem of interest has been termed as well-configuration. Sufficient conditions for a multi-agent system to be well-configured have been provided for different types of directed graphs. For well-configured multi-agent systems, provably correct distributed algorithms have been developed for a number of special cases of the problem. It turns out that the state forms of these algorithms share similarity with so-called matrix-weighted consensus processes \cite{joao,ahn}. Our results imply that the existing sufficient conditions for matrix-weighted consensus, which usually require a tree whose matrix-valued weights are all positive definite, can be significantly relaxed. 

In addition to the two open questions stated at the end of the preceding section, there are a number of directions of future work, including to establish necessary and sufficient conditions for a strongly connected system to be well-configured, to study well-configuration for general rooted graphs, and to derive convergence rates for the proposed algorithms.

\vspace{.1in}

\bibliographystyle{unsrt}
\bibliography{ji}

\end{document}